\input amstex
\magnification=1200
\documentstyle{amsppt}
\NoRunningHeads
\NoBlackBoxes
\topmatter
\title Perception games, the image understanding and interpretational 
geometries
\endtitle
\author Denis V. Juriev
\endauthor
\affil ul.Miklukho-Maklaya 20-180, Moscow 117437 Russia\linebreak
(e-mail: denis\@juriev.msk.ru)
\endaffil
\date math.OC/9905165
\enddate
\keywords Interactive games, Perception, Image Understanding, Dialogues,
Interpretational Geometry\endkeywords
\subjclass 90D20 (Primary) 90D80, 49N55, 93C41, 93B52, 51D05, 68U07 
(Secondary)
\endsubjclass
\abstract\nofrills The interactive game theoretical approach to the
description of perception processes is proposed. The subject is treated 
formally in terms of a new class of the verbalizable interactive games which
are called the {\it perception games}. An application of the previously 
elaborated formalism of dialogues and verbalizable interactive games to the 
visual perception allows to combine the linguistic (such as formal grammars), 
psycholinguistic and (interactive) game theoretical methods for analysis of 
the image understanding by a human that may be also useful for the elaboration 
of computer vision systems. By the way the interactive game theoretical 
aspects of interpretational geometries are clarified.
\endabstract
\endtopmatter
\document

The mathematical formalism of interactive games, which extends one of ordinary
games [1] and is based on the concept of an interactive control, was recently
proposed by the author [2] to take into account the complex composition of 
controls of a real human person, which are often complicated couplings of 
his/her cognitive and known controls with the unknown subconscious behavioral 
reactions. This formalism is applicable also to the description of external 
unknown influences and, thus, is useful for problems in computer science 
(e.g. the semi-artificially controlled distribution of resources) and 
mathematical economics (e.g. the financial games with unknown dynamical 
factors). 

However, the original impetus for the investigations lay in the sphere of
human visual perception [3]. When the first steps were made it became clear
that it is important to understand this sphere formally including it into 
the framework of the elaborating interactive game theory. The interactive
game theoretical definition of dialogues as psycholinguistic phenomena and
the description of the verbalizable interactive games [4] were crucial to
make it possible. This article is an attempt to solve the prescribed problem.
\pagebreak

\head I. Interactive games and their verbalization\endhead

\subhead 1.1. Interactive systems and intention fields\endsubhead

\definition{Definition 1 [2]} An {\it interactive system\/} (with $n$
{\it interactive controls\/}) is a control system with $n$ independent 
controls coupled with unknown or incompletely known feedbacks (the feedbacks
as well as their couplings with controls are of a so complicated nature that 
their can not be described completely). An {\it interactive game\/} is a game 
with interactive controls of each player.
\enddefinition

Below we shall consider only deterministic and differential interactive
systems. In this case the general interactive system may be written in the 
form:
$$\dot\varphi=\Phi(\varphi,u_1,u_2,\ldots,u_n),\tag1$$
where $\varphi$ characterizes the state of the system and $u_i$ are
the interactive controls:
$$u_i(t)=u_i(u_i^\circ(t),\left.[\varphi(\tau)]\right|_{\tau\leqslant t}),$$
i.e. the independent controls $u_i^\circ(t)$ coupled with the feedbacks on
$\left.[\varphi(\tau)]\right|_{\tau\leqslant t}$. One may suppose that the
feedbacks are integrodifferential on $t$.

\proclaim{Proposition [2]} Each interactive system (1) may be transformed
to the form (2) below (which is not, however, unique):
$$\dot\varphi=\tilde\Phi(\varphi,\xi),\tag2$$
where the magnitude $\xi$ (with infinite degrees of freedom as a rule) 
obeys the equation
$$\dot\xi=\Xi(\xi,\varphi,\tilde u_1,\tilde u_2,\ldots,\tilde u_n),\tag3$$
where $\tilde u_i$ are the interactive controls of the form $\tilde 
u_i(t)=\tilde u_i(u_i^\circ(t); \varphi(t),\xi(t))$ (here the dependence
of $\tilde u_i$ on $\xi(t)$ and $\varphi(t)$ is differential on $t$, i.e.
the feedbacks are precisely of the form 
$\tilde u_i(t)=\tilde u_i(u_i^\circ(t);\varphi(t),\xi(t),
\dot\varphi(t),\dot\xi(t),\ddot\varphi(t),\ddot\xi(t),\ldots,
\varphi^{(k)}(t),\mathbreak\xi^{(k)}(t))$).
\endproclaim

\remark{Remark 1} One may exclude $\varphi(t)$ from the feedbacks in
the interactive controls $\tilde u_i(t)$. One may also exclude the
derivatives of $\xi$ and $\varphi$ on $t$ from the feedbacks.
\endremark

\definition{Definition 2 [2]} The magnitude $\xi$ with its dynamical equations
(3) and its cont\-ri\-bution into the interactive controls $\tilde u_i$ will 
be called the {\it intention field}.
\enddefinition

Note that the theorem holds true for the interactive games. In practice, the 
intention fields may be often considered as a field-theoretic description of 
subconscious individual and collective behavioral reactions. However, they 
may be used also the accounting of unknown or incompletely known external 
influences. Therefore, such approach is applicable to problems of computer 
science (e.g. semi-automatically controlled resource distribution) or 
mathematical economics (e.g. financial games with unknown factors).
The interactive games with the differential dependence of feedbacks are
called {\it differential}. Thus, the theorem states a possibility of
a reduction of any interactive game to a differential interactive game
by introduction of additional parameters -- {\sl the intention fields}.

\subhead 1.2. Some generalizations\endsubhead The interactive games introduced 
above may be generalized in the following ways. 

The first way, which leads to the {\it indeterminate interactive games},
is based on the idea that the pure controls $u_i^\circ(t)$ and the 
interactive controls $u_i(t)$ should not be obligatory related in the
considered way. More generally one should only postulate that there are
some time-independent quantities $F_\alpha(u_i(t),u_i^\circ(t),\varphi(t),
\ldots,\varphi^{(k)}(t))$ for the independent magnitudes $u_i(t)$ and 
$u_i^\circ(t)$. Such claim is evidently weaker than one of Def.1. For 
instance, one may consider the inverse dependence of the pure and 
interactive controls: $u_i^\circ(t)=u_i^\circ(u_i(t),\varphi(t),\ldots,
\varphi^{(k)}(t))$.

The second way, which leads to the {\it coalition interactive games}, is
based on the idea to consider the games with coalitions of actions and to
claim that the interactive controls belong to such coalitions. In this case
the evolution equations have the form
$$\dot\varphi=\Phi(\varphi,v_1,\ldots,v_m),$$
where $v_i$ is the interactive control of the $i$-th coalition. If the 
$i$-th coalition is defined by the subset $I_i$ of all players then
$$v_i=v_i(\varphi(t),\ldots,\varphi^{(k)}(t),u^\circ_j| j\in I_i).$$
Certainly, the intersections of different sets $I_i$ may be non-empty so
that any player may belong to several coalitions of actions. Def.1 gives the
particular case when $I_i=\{i\}$.

The coalition interactive games may be an effective tool for an analysis of
the collective decision making in the real coalition games that spread the
applicability of the elaborating interactive game theory to the diverse 
problems of sociology. 

\subhead 1.3. Differential interactive games and their 
$\varepsilon$--representations\endsubhead 

\definition{Definition 3 [4]} The {\it $\varepsilon$--representation\/} of 
differential interactive game is a representation of the differential
feedbacks in the form
$$u_i(t)=u_i(u_i^\circ,\varphi(t),\ldots,\varphi^{(k)}(t);
\varepsilon_i(t))\tag4$$
with the {\sl known\/} function $u_i$ of all its arguments, where
the magnitudes $\varepsilon_i(t)\in\Cal E$ are {\sl unknown\/} functions of
$u_i^\circ$ and $\varphi(t)$ with its higher derivatives:
$$\varepsilon_i(t)=\varepsilon_i(u_i^\circ(t),\varphi(t),\dot\varphi(t),
\ldots,\varphi^{(k)}(t)).$$
\enddefinition

It is interesting to consider several different $\varepsilon$-representations
simultaneously. For such simultaneous $\varepsilon$-representations
with $\varepsilon$-parameters $\varepsilon_i^{(\alpha)}$ a crucial role is
played by the time-independent relations between them:
$$F_\beta(\varepsilon_i^{(1)},\ldots,\varepsilon_i^{(\alpha)},\ldots,
\varepsilon_i^{(N)}; u_i^\circ,\varphi,\ldots,\varphi^{(k)})\equiv0,$$
which are called the {\it correlation integrals}. Certainly, in practice
the correlation integrals are determined {\sl a posteriori\/} and, thus they 
contain an important information on the interactive game. Using the sufficient
number of correlation integrals one is able to construct various algebraic 
structures in analogy to the correlation functions in statistical physics 
and quantum field theory.

\subhead 1.4. Dialogues as interactive games. The verbalization\endsubhead

Dialogues as psycholinguistic phenomena can be formalized in terms of
interactive games. First of all, note that one is able to consider
interactive games of discrete time as well as interactive games of
continuous time above.

\definition{Defintion 4A (the na{\"\i}ve definition of dialogues) [4]}
The {\it dialogue\/} is a 2-person interactive game of discrete time with 
intention fields of continuous time.
\enddefinition

The states and the controls of a dialogue correspond to the speech whereas 
the intention fields describe the understanding. 

Let us give the formal mathematical definition of dialogues now.

\definition{Definition 4B (the formal definition of dialogues) [4]}
The {\it dialogue\/} is a 2-person interactive game of discrete time of 
the form
$$\varphi_n=\Phi(\varphi_{n-1},\vec v_n,\xi(\tau)| 
t_{n-1}\!\leqslant\!\tau\!\leqslant\!t_n).\tag5$$
Here $\varphi_n\!=\!\varphi(t_n)$ are the states of the system at the
moments $t_n$ ($t_0\!<\!t_1\!<\!t_2\!<\!\ldots\!<\!t_n\!<\!\ldots$), 
$\vec v_n\!=\!\vec v(t_n)\!=\!(v_1(t_n),v_2(t_n))$ are the interactive 
controls at the same moments; $\xi(\tau)$ are the intention fields of 
continuous time with evolution equations
$$\dot\xi(t)=\Xi(\xi(t),\vec u(t)),\tag6$$
where $\vec u(t)=(u_1(t),u_2(t))$ are continuous interactive controls with 
$\varepsilon$--represented couplings of feedbacks:
$$u_i(t)=u_i(u_i^\circ(t),\xi(t);\varepsilon_i(t)).$$
The states $\varphi_n$ and the interactive controls $\vec v_n$ are certain
{\sl known\/} functions of the form
$$\aligned
\varphi_n=&\varphi_n(\vec\varepsilon(\tau),\xi(\tau)| 
t_{n-1}\!\leqslant\!\tau\!\leqslant\!t_n),\\
\vec v_n=&\vec v_n(\vec u^\circ(\tau),\xi(\tau)|
t_{n-1}\!\leqslant\!\tau\!\leqslant\!t_n).
\endaligned\tag7
$$
\enddefinition

Note that the most nontrivial part of mathematical formalization of dialogues
is the claim that the states of the dialogue (which describe a speech) are 
certain ``mean values'' of the $\varepsilon$--parameters of the intention
fields (which describe the understanding).

\remark{Important}
The definition of dialogue may be generalized on arbitrary number of players
and below we shall consider any number $n$ of them, e.g. $n=1$ or $n=3$, 
though it slightly contradicts to the common meaning of the word ``dialogue''.
\endremark

An embedding of dialogues into the interactive game theoretical picture
generates the reciprocal problem: how to interpret an arbitrary differential
interactive game as a dialogue. Such interpretation will be called the
{\it verbalization}.

\definition{Definition 5 [4]} A differential interactive game of the form
$$\dot\varphi(t)=\Phi(\varphi(t),\vec u(t))$$
with $\varepsilon$--represented couplings of feedbacks 
$$u_i(t)=u_i(u^\circ_i(t),\varphi(t),\dot\varphi(t),\ddot\varphi(t),\ldots
\varphi^{(k)}(t);\varepsilon_i(t))$$
is called {\it verbalizable\/} if there exist {\sl a posteriori\/}
partition $t_0\!<\!t_1\!<\!t_2\!<\!\ldots\!<\!t_n\!<\!\ldots$ and the 
integrodifferential functionals
$$\aligned
\omega_n&(\vec\varepsilon(\tau),\varphi(\tau)|
t_{n-1}\!\leqslant\!\tau\!\leqslant\!t_n),\\
\vec v_n&(\vec u^\circ(\tau),\varphi(\tau)|
t_{n-1}\!\leqslant\!\tau\!\leqslant\!t_n)
\endaligned\tag8$$ 
such that
$$\omega_n=\Omega(\omega_{n-1},v_n;\varphi(\tau)|
t_{n-1}\!\leqslant\!\tau\!\leqslant\!t_n).
\tag 9$$
\enddefinition

The verbalizable differential interactive games realize a dialogue in sense
of Def.4.

The main heuristic hypothesis is that all differential interactive games
``which appear in practice'' are verbalizable. The verbalization means that 
the states of a differential interactive game are interpreted as intention 
fields of a hidden dialogue and the problem is to describe such dialogue 
completely. If a differential interactive game is verbalizable one 
is able to consider many linguistic (e.g. the formal grammar of a related 
hidden dialogue) or psycholinguistic (e.g. the dynamical correlation of 
various implications) aspects of it.

During the verbalization it is a problem to determine the moments $t_i$. A 
way to the solution lies in the structure of $\varepsilon$-representation.
Let the space $E$ of all admissible values of $\varepsilon$-parameters be
a CW-complex. Then $t_i$ are just the moments of transition of the 
$\varepsilon$-parameters to a new cell. 

\head II. Perception games and the image understanding\endhead

Let us considered a verbalizable interactive game. We shall suppose for 
simplicity that the concrete set is finished if some quantity $F(\omega_n,
\varphi(t))$ reaches some critical value $F_0$. The game will be called 
perception game iff the moments $t_i$ are just the moments of finishing of 
the concrete sets so the multistage perception game realizes a sequence 
of sets with initial states coinciding with the final state of the preceeding
set. Such construction is not senseless contrary to the most of the 
ordinary games because the quantity $F$ should be recalculated with the new 
$\omega$. Thus, we have the following general definition.

\definition{Definition 6} The {\it perception game\/} is a multistage 
verbalizable game (no matter finite or infinite) for which the intervals 
$[t_i,t_{i+1}]$ are just the sets. The conditions of their finishing 
depends only on the current value of $\varphi$ and the state of $\omega$ 
at the beginning of the set. The initial position of the set is the final 
position of the preceeding one.
\enddefinition

Practically, the definition describes the discrete character of the
perception and the image understanding. For example, the goal of a concrete
set may be to perceive or to understand certain detail of the whole image.
Another example is a continuous perception of the moving or changing object.

Note that the definition of perception games is applicable to various forms 
of perception. However, the most interesting one is the visual perception.
Besides the numerous problems of human visual perception of reality 
(as well as of computer vision) there exists a scope of numerous questions 
of the human behaviour in the computer modelled worlds, e.g. constructed 
by use of the 
so-called ``virtual reality'' (VR) technology. There are no an evident 
boundary between them because we can always interpret the internal space 
of our representations as a some sort of the natural ``virtual reality'' 
and apply the analysis of perception in VR to the real image understanding 
as well as to the activity of imagination. So one should convince that it 
is impossible to explain all phenomena of our visual perception of reality 
without deep analysis of its peculiarities in the computer modelled worlds
(cf.[5]). Especially crucial role is played by the so-called {\it integrated
realities\/} (IR), in which the channels only of some kinds of perception are
virtual (e.g. visual) whereas others are real (e.g. tactile, kinesthetic).

The proposed definition allows to take into account the dialogical character
of the image understanding and to consider the visual perception, image
understanding and the verbal (and nonverbal) dialogues together. It may be
extremely useful for the analysis of collective perception, understanding 
and controlling processes in the dynamical environments -- sports, dancings, 
martial arts, the collective controlling of moving objects, etc. 

On the other hand this definition explicates the self-organizing features
of human perception, which may be unraveled by the game theoretical analysis.

And, finally, the definition put a basis for a systematical application of
the linguistic (e.g. formal grammars) and psycholinguistic methods to the 
image understanding as a verbalizable interactive game with a mathematical 
rigor. 

Also interpreting perception processes and the image understanding as the
verbalizable interactive games we obtain an opportunity to adapt some
procedures of the image understanding to the verbalizable interactive
games of a different nature, e.g. to the verbal dialogues. It may enlight
the processes of generation of subjective figurative representations, 
which is important for the analysis of the understanding of speech in
dialogues. 

Traditionally the problems of the visual perception are related to geometry 
(as descriptive as abstract) so it is reasonable to pay an attention to 
the geometrical background for the perception videogames and then to combine 
both geometrical and interactive game theoretical approaches to the visual 
perception and the image understanding.

\head III. Interpretational geometries, intentional anomalous virtual 
realities and their interactive game theoretical aspects\endhead

\subhead 3.1. Interpretational figures [6]\endsubhead
Geometry described below is related to a class of interactive
information systems. Let us call an interactive information system
computer graphic (or interactive information videosystem) if the
information stream ``computer--user'' is organized as a stream of
geometric graphical data on a screen of monitor; an interactive
information system will be called psychoinformation if an information
transmitted by the channel ``user--computer'' is (completely or
partially) subconscious. In general, an investigation of
interactively controlled (psychoinformation) systems for an
experimental and a theoretical explication of possibilities
contained in them, which are interesting for mathematical sciences
themselves, and of ``hidden'' abstract mathematical objects, whose
observation and analysis are actually and potentially realizable
by these possibilities, is an important problem itself. So below
there will be defined the notions of an interpretational figure
and its symbolic drawing that undoubtly play a key role in the
description of a computer--geometric representation of mathematical
data in interactive information systems. Below, however, the accents
will be focused a bit more on applications to informatics preserving
a general experimentally mathematical view, the interpretational
figures (see below) will be used as pointers to droems and interactive
real-time psychoinformation videosystems will be regarded as components
of integrated interactive videocognitive systems for accelerated
nonverbal cognitive communications.

In interactive information systems mathematical data exist in the
form of an interrelation between the geometric internal image (figure)
in the subjective space of the observer and the computer-graphic external
representation. The latter includes visible (drawings of the figure)
and invisible (analytic expressions and algorithms for constructing
these images) elements. Identifying geometric images (figures) in
the internal space of the observer with computer-graphic representations
(visible and invisible elements) is called a {\it translation}, in this
way the visible object may be not identical with the figure, so that
separate visible elements may be considered as modules whose translation
is realized independently. The translation is called an {\it interpretation}
if the translation of separate modules is performed depending on the
results of the translation of preceding ones.

\definition{Definition 7} The figure obtained as a result of
interpretation is called an {\it interpretational figure}.
\enddefinition

Note that the interpretational figure may have no usual formal
definition; namely, only if the process of interpretation admits an
equivalent process of compilation definition of the figure is reduced
to definitions of its drawings that is not true in general. So the
drawing of an interpretational figure defines only dynamical ``technology
of visual perception'' but not its ``image'', such drawings will be
called {\it symbolic}.

The computer-geometric description of mathematical data in interactive
information systems is closely connected with the concept of anomalous
virtual reality.

\subhead 3.2. Intentional anomalous virtual realities [3,6]\endsubhead

\definition{Definition 8}
{\bf (A).} {\it Anomalous virtual reality\/} ({\it AVR\/}) {\it in a narrow
sense\/} means some system of rules of a nonstandard descriptive geometry
adapted for realization on videocomputers (or multisensorial systems of
``virtual reality''). {\it Anomalous virtual reality in a wide
sense\/} also involves an image in cyberspace formed in accordance with
said system of rules. We shall use the term in its narrow sense.
{\bf (B).} {\it Naturalization\/} is the constructing of an AVR from
some abstract geometry or physical model. We say that anomalous virtual
reality {\it naturalizes\/} the abstract model and the model {\it
transcendizes\/} the naturalizing anomalous virtual reality.
{\bf (C).} {\it Visualization\/} is the constructing of certain image
or visual dynamics in some anomalous virtual reality (realized by
hardware and software of a computer-grafic interface of the concrete
videosystem) from the objects of an abstract geometry or processes in
a physical model.
{\bf (D).} Anomalous virtual reality, whose objects depend on the observer,
is called an {\it intentional anomalous virtual reality\/} ({\it IAVR\/}).
The generalized perspective laws for IAVR contain the interactive dynamical
equations for the observed objects in addition to standard (geometric)
perspective laws. In IAVR the observation process consists of a physical
process of observation and a virtual process of intentional governing of
the evolution of images in accordance with the dynamical perspective laws.
\enddefinition

In intentional anomalous virtual reality (IAVR) that is realized
by hardware and software of the {\it computer-graphic interface of
the interactive videosystem} being geometrically modelled by this
IAVR (on the level of descriptive geometry whereas the model transcendizing
this IAVR realizes the same on the level of abstract geometry) respectively,
the observed objects are demonstrated as connected with the observer who
acts on them and determines, or fixes, their observed states so that
the objects are thought only as a potentiality of states from the given
spectrum whose realization depends also on the observer. The symbolic
drawings of interpretational figures may be considered as states of
some IAVR.

Note that mathematical theory of anomalous virtual realities (AVR)
including the basic procedures of naturalization and thanscending
connected AVR with the abstract geometry is a specific branch of {\it
modern nonclassical descriptive (computer) geometry}.

\definition{Definition 8E} The set of all continuously distributed
visual charcteristics of the image in anomalous virtual reality is called
an {\it anomalous color space\/}; the anomalous color space elements
of noncolor nature are called {\it overcolors}, and the quantities
transcendizing them in an abstract model are called {\it ``latent
lights''}. The set of the generalized perspective laws in a fixed
anomalous color space is called a {\it color-perspective system\/}.
\enddefinition

\subhead 3.3. Remarks on the interactive game theoretical aspects 
\endsubhead 
Certainly, the interpretational geometries may be considered
as the perception games. An interesting geometrical consequence of such
approach was proposed [7].

\proclaim{Proposition} {\it There exist models of interpretational
geometries in which there are interpretational figures observed only
in a multi-user mode.}
\endproclaim

It seems that this proposition may be regarded as a startpoint for
the future interactions between geometry and interactive game theory
in the sphere of mathematical foundations for the collective perception 
and image understanding of real objects as well as objects in the computer 
VR or IR systems.

\head IV. Conclusions\endhead

Thus, the interactive game theoretical approach to the description of 
perception processes is proposed. A new class of the multistage verbalizable 
interactive games, the perception games, is introduced. The interactive game 
theoretical aspects of interpretational geometries are clarified. 
Perspectives are sketched.

\Refs
\roster
\item"[1]" Isaaks R., Differential games. Wiley, New York, 1965;\newline
Owen G., Game theory, Saunders, Philadelphia, 1968.
\item"[2]" Juriev D., Interactive games and representation theory. I,II.
E-prints: math.FA/9803020, math.RT/9808098.
\item"[3]" Juriev D., Anomalous color spaces and their structure. Visual
Computer 11(2) (1994) 113-120.
\item"[4]" Juriev D., Interactive games, dialogues and the verbalization.
E-print: math.OC/9903001.
\item"[5]" Forman N., Wilson P., Using of virtual reality for
psychological investigations [in Russian]. Psikhol.Zhurn. 17(2) (1996)
64-79.
\item"[6]" Juriev D., Octonions and binocular mobilevision [in Russian].
Report RCMPI-96/04 (1996); Fund.Prikl.Matem., to appear [Draft English 
e-version: hep-th/9401047];\newline
Juriev D.V., Belavkin--Kolokol'tsov watch--dog effects in interactively 
controlled stochastic dynamical videosystems. Theor.Math.Phys. 106 (1996) 
276-290. Appendix A.
\item"[7]" Juriev D., The advantage of a multi-user mode. E-print:
hep-th/9404137.
\endroster
\endRefs
\enddocument